\documentclass[suppldata]{interact}

\usepackage{epstopdf}
\usepackage[caption=false]{subfig}
\usepackage[doublespacing]{setspace}

\usepackage[numbers,sort&compress]{natbib}
\bibpunct[, ]{[}{]}{,}{n}{,}{,}
\makeatletter
\def\NAT@def@citea{\def\@citea{\NAT@separator}}
\makeatother

\theoremstyle{plain}

\usepackage{amsmath,amssymb,amsthm}
\usepackage[hyphens]{url} 

\usepackage{multirow}
\usepackage{array}
\usepackage{epsfig,rotating}
\usepackage{latexsym}
\usepackage{bm}
\usepackage{upgreek}
\usepackage{float}
\usepackage{times}

\def\eop{\hfill $\Box$ \medskip}

\newtheorem{theorem}{{Theorem}}
\newtheorem{prop}{Proposition}[section]

\begin{document}

\title{Empirical Likelihood Method for Complete Independence Test on High Dimensional Data}

\author{
\name{Yongcheng Qi\textsuperscript{1} and Yingchao Zhou\textsuperscript{1}}
\affil{\textsuperscript{1}Department of Mathematics and Statistics, University of Minnesota Duluth, 1117 University Drive, Duluth, MN 55812, USA
}
}

\maketitle

\footnotetext[2]{This is an Accepted Manuscript version of the
following article, accepted for publication in {\it Journal of Statistical Computation and Simulation } [https://doi.org/10.1080/00949655.2022.2029860].
 It is deposited under the terms of the Creative
Commons Attribution-NonCommercial License
(http://creativecommons.org/licenses/by-nc/4.0/), which permits
non-commercial re-use, distribution, and reproduction in any medium,
provided the original work is properly cited.}

\begin{abstract}
Given a random sample of size $n$ from a
$p$ dimensional random vector, where both $n$ and $p$ are large, we
are interested in testing whether the $p$ components of the random
vector are mutually independent. This is the so-called complete
independence test.  In the multivariate normal case, it is
equivalent to testing whether the correlation matrix is an identity
matrix. In this paper, we propose a one-sided empirical likelihood method
for the complete independence test for multivariate normal data
based on squared sample correlation coefficients. The limiting
distribution for our one-sided empirical likelihood test statistic
is proved to be $Z^2I(Z>0)$ when both $n$ and $p$ tend to infinity,
where $Z$ is a standard normal random variable. In order to improve the power of the empirical likelihood test statistic, we also introduce a rescaled empirical likelihood test statistic.
We carry out an extensive simulation study to compare the performance of the
rescaled empirical likelihood method and two other statistics which are
related to the sum of squared sample correlation coefficients.

\end{abstract}

\begin{keywords} Empirical likelihood; complete
independence test; high dimension; multivariate normal distribution
\end{keywords}


\newpage

\section{Introduction}\label{introduction}

Statistical inference on high dimensional data
 has gained a wide range of applications in
recent years. New techniques generate a vast collection of data sets
with high dimensions, for example, trading data from financial
market, social network data and biological data like microarray and
DNA data. The dimension of these types of data is not small compared
with sample size, and typically of the same order as sample size or
even larger. Yet classical multivariate statistics usually deal with
data from normal distributions with a large sample size $n$ and a
fixed dimension $p$, and one can easily find some classic treatments
in reference books such as Anderson~\cite{Anderson},
Morrison~\cite{Morrison} and Muirhead~\cite{Muirhead1982}.

Under multivariate normality settings,  the likelihood ratio test
statistic converges in distribution to a chi-squared distribution
when $p$ is fixed.  However, when $p$ changes with $n$ and tends to
infinity, this conclusion is no longer true as discovered in Bai et
al.~\cite{Bai}, Jiang et al.~\cite{Jiang2012},  Jiang and
Yang~\cite{Jiang2013}, Jiang and Qi~\cite{Jiang2015a}, Qi et
al.~\cite{QWZ19}, among others.  The results in these papers
indicate that the chi-square approximation fails when $p$ diverges as
$n$ goes to infinity.

The test of complete independence of a random vector is to test
whether all the components of the random vector are mutually
independent. In the multivariate normal case, the test of complete
independence is equivalent to the test whether covariance matrix is
a diagonal matrix, or whether the correlation matrix is the identity
matrix.

For more details, we assume $X=(X_1, \cdots, X_p)$ is a random
vector from a $p$-dimensional multivariate normal distribution
$N_p(\bm\mu,\mathbf\Sigma)$, where $\bm\mu$ denotes the mean vector,
and $\mathbf\Sigma$ is a $p\times p$ covariance matrix. Given a
random sample of size $n$ from the normal distribution,
$\mathbf{x}_1, \mathbf{x}_2,\cdots,\mathbf{x}_n$, where
$\mathbf{x}_i=(x_{i1},x_{i2},\cdots,x_{ip})$ for $1\le i\le n$,
Pearson's correlation coefficient between the $i$-th and $j$-th
components is given by
\begin{equation}\label{rij}
r_{ij}=\frac{\sum\limits_{k=1}^n(x_{ki}-\bar{x_i})(x_{kj}-\bar{x}_j)}{\sqrt{\sum\limits_{k=1}^n{(x_{ki}-\bar{x}_i)^2}
\cdot {\sum\limits_{k=1}^n(x_{kj}-\bar{x_j})^2}}},
\end{equation}
where $\bar x_i=\frac{1}{n}\sum\limits_{k=1}^n x_{ki}$ and $\bar
x_j=\frac{1}{n}\sum\limits_{k=1}^n x_{kj}$ for $1\le i, j\le p$. Now
we set $\mathbf{R}_n=(r_{ij})_{p\times p}$ as the sample correlation
coefficient matrix.

The complete independence test for the normal random vector is
\begin{equation}\label{HH}
H_0:\boldsymbol\Gamma=\mathbf{I_p}~~~ vs ~~~H_a:
\boldsymbol\Gamma\neq\mathbf{I}_p,
\end{equation}
where $\boldsymbol\Gamma$ is the population correlation matrix and
$\mathbf{I}_p$ is $p\times p$ identity matrix. When $p<n$, the
likelihood ratio test statistic for \eqref{HH} is a function of
$|\mathbf{R}_n|$, the determinant of $\mathbf{R}_n$, from
Bartlett~\cite{Bartlett1954} or Morrison~\cite{Morrison}.  In
traditional multivariate analysis, when $p$ is a fixed integer, we
have under the null hypothesis in \eqref{HH} that
\[
-(n-1-\frac{2p+5}{6})\log|\mathbf{R}_n|\overset{d}\to\chi^2_{p(p-1)/2}~~~~\mbox{
as }n\to\infty,
\]
where $\chi^2_{f}$ denotes a chi-square distribution with $f$ degrees of freedom.

When $p=p_n$ depends on $n$ with $2\le p_n<n$ and $p_n\to\infty$,
the likelihood ratio method can still be applied to test \eqref{HH}.
The limiting distributions of the likelihood ratio test statistics
in this case have been discussed in the aforementioned papers. It is
worth mentioning that Qi et al.~\cite{QWZ19} propose an adjusted
likelihood ratio test statistic and show that the distribution of
the adjusted likelihood test statistic can be well approximated by a
chi-squared distribution whose number of degrees of freedom depends
on $p$ regardless of whether $p$ is fixed or divergent.

The limitation of the likelihood ratio test is that the dimension
$p$ of the data must be smaller than the sample size $n$. Many other
likelihood tests related to the sample covariance matrix or sample
correlation matrix have the same problem as the sample covariance
matrices are degenerate when $p\ge n$. In order to relax this
constraint, a new test statistic using the sum of squared sample
correction coefficients is proposed by Schott~\cite{Schott2005} as
follows
\[
t_{np}=\sum_{1\le j<i\le p}r_{ij}^2.
\]
Assume that the null hypothesis of \eqref{HH} holds. Under
assumption $\lim\limits_{n\to\infty}p_n/n=\gamma\in (0, \infty)$,
Schott~\cite{Schott2005} proves that $t_{np}-\frac{p(p-1)}{2(n-1)}$
converges in distribution to a normal distribution with mean $0$ and
variance $\gamma^2$, that is,
\begin{equation}\label{schott1}
t_{np}^*:=\frac{t_{np}-\frac{p(p-1)}{2(n-1)}}{\sigma_{np}}\overset{d}{\to}
N(0,1),
\end{equation}
where $\sigma_{np}^2=\frac{p(p-1)(n-2)}{(n-1)^2(n+1)}$.

Recently, Mao~\cite{Mao2014} proposes a different test for complete
independence. His test statistic is closely related to Schott's test
and is defined by
\[
T_{np}=\sum_{1\le j<i\le p}\frac{r_{ij}^2}{1-r_{ij}^2}.
\]
It has been proved in Mao~\cite{Mao2014} that $T_{np}$ is
asymptotically normal under the null hypothesis of \eqref{HH} and the
assumption that $\lim\limits_{n\to\infty}p_n/n=\gamma\in (0, \infty)$.

Very recently, Chang and Qi~\cite{CQ18} investigate the limiting
distributions for the two test statistics above under less
restrictive conditions on $n$ and $p$.   Chang and Qi~\cite{CQ18}
show that \eqref{schott1} is also valid under the general condition
that $p_n\to\infty$ as $n\to\infty$, regardless of the convergence
rate of $p_n$.  Thus, the normal approximation in \eqref{schott1}
based  on $t_{np}^*$ yields an approximate level $\alpha$ rejection
region
\begin{equation}\label{RT*}
\mathcal{R}_t^*(\alpha)=\Big\{t_{np}\ge
\frac{p(p-1)}{2(n-1)}+z_{1-\alpha}\sqrt{\frac{p(p-1)(n-1)}{(n-1)^2(n+1)}}\Big\},
\end{equation}
where $z_{\alpha}$ is a $\alpha$ level critical value of the
standard normal distribution.

Furthermore, Chang and Qi~\cite{CQ18} propose adjusted test
statistics whose distribution can be fitted by chi-squared
distribution regardless of how $p$ changes with $n$ as long as $n$
is large. Chang and Qi's~\cite{CQ18}  adjusted test statistics
$t_{np}^c$ is defined as
\begin{equation}\label{tnc}
t_{np}^c=\sqrt{p(p-1)}t_{np}^*+\frac{p(p-1)}{2}.
\end{equation}
Chang and Qi show that
\[
\sup\limits_x\left|P(t_{np}^c \leqslant
x)-P(\chi^2_{p(p-1)/2}\leqslant x)\right|\rightarrow 0
\]
as long as $p_n\to\infty$ as $n\to\infty$. 
Let
$\chi^2_f(\alpha)$ denote the $\alpha$ level critical value of
$\chi^2_f$. Then an approximate level $\alpha$ rejection region
based on $t_{np}^c$ is given by
\begin{equation}\label{RTC}
\mathcal{R}_t^c(\alpha)=\Big\{t_{np}\ge
\frac{p(p-1)}{2}(1-\sqrt{\frac{n-2}{n+1}})+
\chi^2_{\frac{p(p-1)}{2}}(\alpha)\sqrt{\frac{n-2}{(n-1)^2(n+1)}}\Big\}.
\end{equation}

One can find more references on test for complete independence in
Mao~\cite{Mao2014} or Chang and Qi~\cite{CQ18}.

In practice,  the assumption of normality for distributions may be violated.
Now we assume $X=(X_1, \cdots, X_p)$ is a random
vector and $X_1, \cdots, X_p$ are identically distributed with distribution function $F$.
Given a random sample of size $n$, $\mathbf{x}_1, \mathbf{x}_2,\cdots,\mathbf{x}_n$, where
$\mathbf{x}_i=(x_{i1},x_{i2},\cdots,x_{ip})$ for $1\le i\le n$, are drawn from the distribution of $X=(X_1, \cdots, X_p)$,  and define
Pearson's correlation coefficients $r_{ij}$'s as in \eqref{rij}. By using the Stein method,  Chen and Shao~\cite{CS2012} show that
\eqref{schott1} holds under some moment conditions of $F$ if $p_n/n$ is bounded.

In this paper, we propose to apply empirical likelihood method to
the testing problem \eqref{HH}. The empirical likelihood is a
nonparametric statistical method proposed by Owen~\cite{Owen1988,
Owen1990}, which is originally used to test the mean vector of a
population based on a set of independent and identically distributed
(i.i.d.) random variables. Empirical likelihood does not require to
specify the family of distributions for the data and it possesses
some good properties of the likelihood methods.


The rest of the paper is organized as follows. In
Section~\ref{main}, we first introduce a one-sided empirical
likelihood method for the mean of a set of random variables with a
common mean and then establish the connection between the test of
complete independence and the one-sided empirical likelihood method.
Our main result concerning the limiting distribution of the
one-sided empirical likelihood ratio statistic is also given in
Section~\ref{main}.  In Section~\ref{sim}, we carry out a simulation
study to compare the performance of the empirical likelihood method
and normal approximation based on Schott's test statistic and
chi-square approximation based on Chang and Qi's adjusted test
statistic.  In our simulation study, we also apply these methods to some other distributions such as the exponential
distributions and mixture of the exponential and normal
distributions so as to compare their adaptability to non-normality.
The proofs of the main results are given in Section~\ref{proofs}.

\section{Main Results}\label{main}

In this section, we apply the empirical likelihood method to the test of complete independence.
First, we assume $X=(X_1, \cdots, X_p)$ is a random
vector from a $p$-dimensional multivariate normal distribution.
Under the null hypothesis of
\eqref{HH}, $\{r_{ij}^2, 1\le i<j\le p\}$ are random variables from an identical distribution with mean $\frac{1}{n-1}$.
As a matter of fact, it follows from Corollary~5.1.2 in Muirhead~\cite{Muirhead1982} that $r_{ij}^2$ has the same distribution as $T^2/(n-2+T^2)$ under the null hypothesis of
\eqref{HH}, where $T$ is a random variable having $t$-distribution with $n-2$ degrees of freedom.
 $\{r_{ij}^2, 1\le i<j\le p\}$ are asymptotically independent if the sample size $n$ is large.
We will
develop a one-sided empirical likelihood test statistic and apply it
to the data set $\{(n-1)r_{ij}^2, 1\le i<j\le p\}$, where $p=p_n$ is a
sequence of positive integers such that $p_n\to\infty$ as
$n\to\infty$.  As an extension, we then consider the case when $X=(X_1, \cdots, X_p)$ is a random
vector with an identical marginal distribution function $F$ which is not necessarily Gaussian. When the $p$ components of $X$ are independent, we demonstrate that the empirical likelihood method we develop under normality works for general distribution $F$ as well if some additional conditions are satisfied.

\subsection{One-sided empirical likelihood test}

Consider a random sample of size $N$, namely $y_1,\cdots, y_N$.
Assume the sample comes from a population with mean $\mu$ and
variance $\sigma^2$. The empirical likelihood function for the mean
$\mu$ is defined as
\begin{equation}\label{L}
L(\mu)=\sup\Big\{\prod^N_{i=1}\omega_i\Big|
\sum^N_{i=1}\omega_iy_i=\mu,  \omega_i\ge 0,
\sum^N_{i=1}\omega_i=1\Big\}.
\end{equation}
The function $L(\mu)$ is well defined if $\mu$ belongs to the convex hull given by
\[
H:=\Big\{\sum^N_{i=1}\omega_iy_i\Big|\sum^N_{i=1}\omega_i=1, \omega_i>0, ~i=1, \cdots, N\Big\};
\]
otherwise, set $L(\mu)=0$.   We see that $H=(\min\limits_{1\le i\le N} y_i, \max\limits_{1\le i\le N} y_i)$.

Assume $\mu\in H$.
By the standard Lagrange multiplier technique,  the supremum on the right-hand side of \eqref{L} is achieved at
\begin{equation}\label{omega-i}
\omega_i=\frac{1}{N(1+\lambda(y_i-\mu))},~~~ i=1, \cdots, N,
\end{equation}
where $\lambda$ is the solution to  equation $g(\lambda)=0$, with $g(\lambda)$ defined as follows
 \begin{equation}\label{lambda}
g(\lambda):=\sum_{i=1}^{N}\frac{y_i-\mu}{1+\lambda(y_i-\mu)}.
\end{equation}
Assume $\min\limits_{1\le i\le N} y_i< \max\limits_{1\le i\le N} y_i$.  When
$\mu\in H$, then the function $g(\lambda)$ defined in \eqref{lambda}
is strictly increasing for $\displaystyle\lambda\in (-(\max_{1\le
i\le N} y_n-\mu)^{-1}, (\mu-\min_{1\le i\le N} y_i)^{-1})$. A
solution to $g(\lambda)=0$ in this range exists and the solution must
be unique.

\begin{prop}\label{prop1} Assume $y_1,\cdots, y_N$ are $N$ observations with $y_i\ne y_j$
for some $i$ and $j$. Then $\log L(\mu)$ is strictly concave in $H$, and
$L(\bar{y})=\sup\limits_{\mu} L(\mu)=N^{-N}$,  where $\bar{y}=\frac1N\sum\limits_{i=1}^N y_i$.
\end{prop}

\vspace{10pt}

\noindent \textbf{Remark}.  The results in Proposition~\ref{prop1}
are well-known among the researchers in the area of empirical
likelihood methods. A short proof will be given in
Section~\ref{proofs} for completeness.

Consider the following two-sided test problem
\[
H_0:\mu=\mu_0~~~vs~~~~H_a:\mu\ne\mu_0.
\]
The empirical likelihood ratio is given by
\[
\frac{L(\mu_0)}{\sup\limits_{\mu\in R}L(\mu)}=\frac{L(\mu_0)}{N^{-N}}=\prod^N_{i=1}(1+\lambda(y_i-\mu_0))^{-1},
\]
where $\lambda$ is the solution to the following equation
 \[
\frac{1}{N}\sum_{i=1}^{N}\frac{y_i-\mu_0}{1+\lambda(y_i-\mu_0)}=0.
\]
Therefore, the log-empirical likelihood test statistic is given by
\begin{equation}\label{lmu_0}
\ell(\mu_0):=-2\log\frac{L(\mu_0)}{\sup\limits_{\mu\in
R}L(\mu)}=2\sum^N_{i=1}\log(1+\lambda(y_i-\mu_0)).
\end{equation}

It is proved in Owen~\cite{Owen2001} that $\ell(\mu_0)$ converges in
distribution to a chi-square distribution with one degree of freedom
if $y_1,\cdots, y_N$ are i.i.d. random variables with mean $\mu_0$
and a finite second moment.

 Our interest here is to consider a
one-sided test
\begin{equation}\label{one-sided}
H_0:\mu=\mu_0~~~vs~~~~H_a:\mu>\mu_0.
\end{equation}
According to Proposition~\ref{prop1},  $L(\mu)$ is increasing in $(-\infty,\bar{y})$ and decreasing in $(\bar{y}, \infty)$, which implies
$\sup\limits_{\mu\ge \mu_0}L(\mu)=L(\mu_0)I(\bar{y}<\mu_0)+N^{-N}I(\bar{y}\ge \mu_0)$. Therefore, the empirical likelihood ratio corresponding to test \eqref{one-sided} is
\[
\frac{L(\mu_0)}{\sup\limits_{\mu\ge \mu_0}L(\mu)}
=\left\{
\begin{array}{ll}
   \frac{L(\mu_0)}{N^{-N}}, & \hbox{ if } \bar{y}\ge \mu_0;\\
   1, & \hbox{ if } \bar{y}<\mu_0.
   \end{array}
   \right.
\]
Then the log-empirical likelihood test statistic for test \eqref{one-sided} is
\begin{equation}\label{one-sided-L}
\ell_n(\mu_0):=-2\log \frac{L(\mu_0)}{\sup\limits_{\mu\ge
\mu_0}L(\mu)}=\ell(\mu_0)I(\bar y\ge \mu_0),
\end{equation}
where $\ell(\mu_0)$ is defined in \eqref{lmu_0}.

\subsection{Empirical likelihood method for testing complete independence}

Let $r$ denote the sample Pearson correlation coefficient based on a
random sample of size $n$ from a bivariate normal distribution with
correlation coefficient $\rho$.  From Muirhead~\cite{Muirhead1982},
page 156,
 \[
E(r^2)=1-\frac{n-2}{n-1}(1-\rho^2)_{2}F_1(1,1;\frac{1}{2}n+1;\rho^2),
 \]
 where
 \[
 _{2}F_1(a,b;c;z)=1+\frac{ab}{1!c}z+\frac{a(a+1)b(b+1)}{2!c(c+1)}z^2+\dots=1+\sum\limits_{k=1}^{\infty}\frac{(a)_k(b)_k}{(c)_k}\frac{z^k}{k!}\]
 is the hypergeometric function,  $(a)_k=\Gamma(a+k)/\Gamma(a)$, and
$\Gamma(x)=\int^\infty_0t^{x-1}e^{-t}dt$ is the gamma function. It is easy to check when $\rho=0$, $_{2}F_1(1,1;\frac{1}{2}n+1;\rho^2)=1$,
and $E(r^2)=1-\frac{n-2}{n-1}=\frac{1}{n-1}$; when $\rho\neq 0$,
$_{2}F_1(1,1;\frac{1}{2}n+1;\rho^2)<
1+\sum\limits_{k=1}^{\infty}\rho^{2k}=\frac{1}{1-\rho^2}$, and thus,
$E(r^2)> 1-\frac{n-2}{n-1}=\frac{1}{n-1}$.

First,  we assume $X=(X_1, \cdots, X_p)$ is a random
vector from a $p$-dimensional multivariate normal distribution
$N_p(\bm\mu,\mathbf\Sigma)$.  Review the sample correlation coefficients $r_{ij}$ defined in
\eqref{rij}.  Denote the correlation matrix of $\bm\Sigma$ by
$\boldsymbol\Gamma=(\gamma_{ij})$. From the above discussion, we
have that under the null hypothesis of \eqref{HH},
$E(r_{ij}^2)=\frac{1}{n-1}$ for all $1\le i<j\le p$,
$E(r_{ij}^2)\ge\frac{1}{n-1}$ under the alternative of \eqref{HH}
and at least one of the inequalities is strict.
We see that
test \eqref{HH} is equivalent to the following one-tailed test
\[
H_0:E(\bar{r}_{ij})=1,~1\le i< j\le p~~vs~~H_a: E(\bar{r}_{ij})>1~\mbox{ for some }1\le 1<j\le p,
\]
where $\bar{r}_{ij}=(n-1)r_{ij}^2$.
Under the null hypothesis of \eqref{HH},  $\{(n-1)r_{ij}^2,~ 1\le i< j\le p\}$ are identically distributed with mean $1$ and
variance $\frac{2(n-2)}{(n+1)}$.
 We also notice from Chang and Qi~\cite{CQ18}
that $\{(n-1)r_{ij}^2,~1\le i< j\le p\}$ behave as if they were independent and
identically distributed. For these reasons, we propose a one-sided
empirical likelihood ratio test as follows.

Rewrite $\{(n-1)r_{ij}^2,~1\le i < j\le p\}$ as $y_1, \cdots, y_N$, where
$N=p(p-1)/2$. Then $y_1,\cdots, y_N$ are asymptotically i.i.d
with mean $1$. Define the one-sided log-empirical
likelihood ratio test statistics as in \eqref{one-sided-L} with
$\mu_0=1$, or equivalently
\begin{equation}\label{l-one}
\ell_n:=\ell_n(1)=2I(\bar r\ge 1)\sum_{1\le
i<j\le p}\log \Big(1+\lambda\big((n-1)r_{ij}^2-1\big)\Big),
\end{equation}
where $\lambda$ is the solution to the equation
\[
\sum_{1\le i<j\le p}\frac{(n-1)r_{ij}^2-1}{1+\lambda\big((n-1)r_{ij}^2-1\big)}=0,
\]
and $\bar r=\bar y=\frac{n-1}{N}\sum_{1\le i<j\le p}r_{ij}^2$.

Our first result on empirical likelihood method for testing the complete independence under normality in the paper is as follows.

\begin{theorem}\label{thm1}
Assume $p=p_n\to\infty$ as $n\to\infty$. Then
$\ell_n\overset{d}{\rightarrow} Z^2I(Z>0)$ as $n\rightarrow \infty$
under the null hypothesis of \eqref{HH}, where $Z$ is a standard
normal random variable.
\end{theorem}

Let $\Phi$ denote the cumulative distribution function of the standard normal distribution, i.e,
\[
\Phi(x)=\frac{1}{\sqrt{2\pi}}\int^x_{-\infty}e^{-t^2/2}dt,  ~~~~x\in (-\infty, \infty).
\]
Let $G$ denote the cumulative distribution function of $Z^2I(Z>0)$. Then
\[
G(x)=\left\{
       \begin{array}{ll}
         0, & \hbox{$x<0$;} \\
        \Phi(\sqrt{x}), & \hbox{$x\ge 0$.}
       \end{array}
     \right.
\]
Therefore, for any $\alpha\in (0,\frac12)$, an $\alpha$-level
critical value of $G$ is given by $z_\alpha^2$, where $z_{\alpha}$
is an $\alpha$-level critical value for the standard normal
distribution.  Based on Theorem~\ref{thm1},  a level $\alpha$
rejection region for test on \eqref{one-sided} is
\begin{equation}\label{Re}
\mathcal{R}_{e}(\alpha)=\Big\{\ell_n\ge z_{\alpha}^2\Big\}.
\end{equation}
Here we only consider $\alpha<\frac{1}{2}$ because $Z^2I(Z>0)$ is
nonnegative, $P(Z^2I(Z>0)>c|H_0)<\frac{1}{2}$ if $c>0$, and
$P(Z^2I(Z>0)>c|H_0)=1$ if $c\le0$.

Now we consider the general case
when $X=(X_1, \cdots, X_p)$ is a random
vector with independent and identically distributed components.
The one-sided empirical likelihood test statistic $\ell_n$ based on $\{(n-1)r_{ij}^2, ~1\le i<j\le p\}$ is defined as in \eqref{l-one}. The limiting distribution for $\ell_n$ is the same as that under normality.

\begin{theorem}\label{thm2}
Assume $X_1, \cdots, X_p$ are independent and identically distributed and  $E(X_1^{24})<\infty$.
If $p=p_n\to\infty$ as $n\to\infty$ and $p_n/n$ is bounded, then
$\ell_n\overset{d}{\rightarrow} Z^2I(Z>0)$ as $n\rightarrow \infty$.
\end{theorem}

Compared with Theorem~\ref{thm1},  $p_n$ in Theorem~\ref{thm2} is restricted in a smaller range and it can be of the same order as $n$.

To demonstrate the performance of empirical likelihood method and two other test statistics, we have a numerical study.  Our simulation study indicates that the empirical likelihood test \eqref{Re} based om $\ell_n$ maintains a very stable size or type I error. In terms of size, $\ell_n$ is more accurate $t_{np}^c$ and $t_{np}^*$. Most of the time, $t_{np}^c$ and $t_{np}^*$ have slightly larger sizes than $0.05$ when the nominal level $\alpha$ is $0.05$, and their powers are also slightly larger than that of $\ell_n$ in our simulation study. For simplicity purpose, the simulation result on $\ell_n$ is not shown in this paper.

In order to balance the size and power for the empirical likelihood method, We introduce a rescaled empirical likelihood statistic, $\bar{\ell}_n$, defined as follows
\begin{equation}\label{bar-ell}
\bar{\ell}_n=\frac{2(n-1)(n+1)}{3(p-1)(p+4)}\ell_n\sum_{1\le i<j\le p}r_{ij}^4.
\end{equation}
Under conditions of Theorems~\ref{thm1} or \ref{thm2}, $\bar{\ell}_n$ and $\ell_n$ have the same limiting distribution, that is,
\begin{equation}\label{bar-ell-limit}
\bar{\ell}_n\overset{d}{\rightarrow} Z^2I(Z>0)~~\mbox{ as  }~~n\rightarrow \infty
\end{equation}
provided that
\begin{equation}\label{1-in-probability}
 \frac{2(n-1)(n+1)}{3(p-1)(p+4)}\sum_{1\le i<j\le p}r_{ij}^4\overset{p}\to 1.
\end{equation}
This equation will be verified in Section~\ref{proofs}. Based on \eqref{bar-ell-limit}, a level $\alpha$ test rejects the complete independence if $\bar{\ell}_n$ falls into the rejection region
\begin{equation}\label{bar-Re}
\bar{\mathcal{R}}_{e}(\alpha)=\Big\{\bar{\ell}_n\ge z_{\alpha}^2\Big\}.
\end{equation}

\section{Simulation}\label{sim}

In this section, we will consider the following three test
statistics for testing complete independence \eqref{HH}, including
Schott's test statistic $t_{np}^*$ given in \eqref{schott1}, Chang
and Qi's adjusted test statistic $t_{np}^c$ defined in \eqref{tnc},
and the rescaled empirical likelihood test statistic $\bar{\ell}_n$ given in
\eqref{bar-ell}. The corresponding rejection regions are given in
\eqref{RT*}, \eqref{RTC}, and \eqref{bar-Re}, respectively.  All
simulations are implemented by the software $\texttt{R}$.

For sample size $n=20,50,100$ and dimension $p=10,20,50,100$, we
apply the three test statistics to each of five distributions for
10000 iterations to obtain the empirical sizes and the empirical
powers of the tests. We set the nominal type I error $\alpha=0.05$.
The five distributions include the normal,  the uniform over $[-1,
1]$, the exponential, the mixture of the normal and exponential
distributions, and the sum of normal and exponential distributions.

To control the dependence structure, we introduce a covariance
matrix $\boldsymbol\Gamma_\rho$ defined by
\begin{equation}\label{Sigma}
\boldsymbol\Gamma_\rho=\big(\gamma_{ij}\big)_{p\times p}, \mbox{
with }\gamma_{ii}=1, \mbox{ and }\gamma_{ij}=\rho \mbox{ if }i\ne j,
\end{equation}
which is also a correlation matrix. In our simulation study,  we
generate random samples from the distribution of a random vector
$X=(X_1, \cdots, X_p)$ with covariance matrix
$\boldsymbol\Gamma_\rho$ or correlation matrix
$\boldsymbol\Gamma_\rho$.
For details, see the five distributions described
below. For all distributions we consider, the
observations have independent components when $\rho=0$ and
positively dependent components when $\rho>0$.  We choose very small
values for $\rho$ such as $\rho=0.02$ and $0.05$. When the
value of $\rho$ is large, the resulting powers for all three methods
will be too close to $1$, and the comparison is meaningless.
Therefore, based on $10,000$ replicates, the sizes for three test
statistics are estimated when $\rho=0$, and their powers are
estimated when $\rho=0.02$ and $0.05$.   All results are
reported in Tables~\ref{table-norm} to \ref{table-sum}.

\vspace{10pt}
\noindent\textbf{a. ~ Normal Distribution}

The observations are drawn from a multivariate normal random vector
$X=(X_1,\cdots, X_p)$ with mean $\boldsymbol{\mu}=(0,\cdots, 0)$ and
variance matrix $\boldsymbol\Gamma_{\rho}$ specified in
\eqref{Sigma}.  The results on the empirical sizes and powers are
given in Table~\ref{table-norm}.

\vspace{10pt}
\noindent\textbf{b. ~Uniform Distribution}

We first generate $p+1$ i.i.d. random variables $Y_0, Y_1, \cdots,
Y_p$ from Uniform $(-1,1)$ distribution, then set
$X_i=\frac{\sqrt{\rho}}{\sqrt{1-\rho}}Y_0+Y_i$, $i=1, 2,\cdots, p$.
It is easy to verify that random vector $X=(X_1,\cdots, X_p)$ has
mean $\boldsymbol{\mu}=(0,\cdots,0)$ and correlation matrix
$\boldsymbol\Gamma_{\rho}$ as defined in \eqref{Sigma}.  The results
on the empirical sizes and powers are given in Table~\ref{table-unif}.

\vspace{10pt}
\noindent\textbf{c. ~Exponential Distribution}

We generate $p+1$ i.i.d. random variables $Y_0, Y_1, \cdots, Y_p$
from the unit exponential distribution, then define
$X_i=\frac{\sqrt{\rho}}{\sqrt{1-\rho}}Y_0+Y_i$, $i=1, 2,\cdots, p$.
The random vector $X=(X_1,\cdots, X_p)$ has a correlation matrix
$\boldsymbol\Gamma_{\rho}$ as defined in \eqref{Sigma} for $\rho\in
[0,1)$.   The results on the empirical sizes and powers are given in
Table~\ref{table-exp}.

\vspace{10pt}
\noindent\textbf{d. ~Mixture of Normal and Exponential Distributions}

The random vector $X=(X_1, \cdots, X_p)$ is sampled from a mixture
of the normal and exponential distributions which is with 90\%
probability from the multivariate normal with mean
$\boldsymbol\mu=(1,\cdots,1)$ and covariance matrix
$\boldsymbol\Gamma_\rho$ given in \eqref{Sigma} and with 10\%
probability from a random vector $(Y_1, \cdots, Y_p)$ where $Y_1,
\cdots, Y_p$ are i.i.d. unit exponential random variables.
 The results on the empirical sizes and powers are given
in Table~\ref{table-mix}.

\vspace{10pt}
\noindent\textbf{e. ~Sum of Normal and Exponential Distribution}

The random vector $X=(X_1,\cdots, X_p)$ is a weighted sum of two
independent random vectors, $U$ and $V$,  $X=U+0.01V$, where $U$ is
from a multivariate normal distribution with mean
$\boldsymbol{\mu}=(0,\cdots,0)$ and covariance matrix
$\boldsymbol\Gamma_{\rho}$ defined in \eqref{Sigma}, and $V=(Y_1,
\cdots, Y_p)$ with $Y_i$'s being i.i.d. unit exponential random
variables.  The results on the empirical sizes and powers are given
in Table~\ref{table-sum}.

\vspace{10pt}

From the simulation results, the empirical sizes for all three tests are close to $0.05$ which is the nominal type I error we set in the simulation, especially when both $n$ and $p$ are large. Test statistic $t_{np}^c$ has the smallest size in most cases, and it is a little bit conservative sometimes.
The size of $\bar{\ell}_n$ is between that of $t_{np}^c$ and $t_{np}^*$ and both $\bar{\ell}_n$ and $t_{np}^*$ are comparable
for most combinations of $n$ and $p$.

As we expect, the powers of all three test statistics become higher
as $p$ grows larger. The increase in $n$ also brings about an
increase in power, but not as much as the increase in $p$ does,
because $\frac{p(p-1)}{2}$ is the number of $r_{ij}^2$'s involved in
the test. All test statistics achieve high power when $\rho=0.05$.
Three test statistics result in comparable powers in general,
although the power of Chang and
Qi's test statistic is occasionally
a little bit less than the other two test statistics. These
differences may be due to the fact that Chang and
Qi's test statistic maintain a lower type I error.

In summary, in this paper, we have developed the one-sided empirical likelihood method and proposed the rescaled empirical likelihood test statistic for testing the complete independence for high dimensional random vectors. The rescaled empirical likelihood test statistic performs very well in terms of the size and power and can serve as a good alternative to the existent test statistics in the literature.

\section{Proofs}\label{proofs}

\noindent{\it Proof of Proposition~\ref{prop1}.} To prove the strict
concavity of $L(\mu)$, we need to show that for $\mu_1, \mu_2\in H$,
$\mu_1\ne \mu_2,$
\begin{equation}\label{concave}
\log L(t\mu_1+(1-t)\mu_2)>t\log L(\mu_1)+(1-t)\log L(\mu_2),~~ t\in
(0,1).
\end{equation}
Since $\mu_j\in H$ for $j=1,2$,  we have $\log
L(\mu_j)=\log\prod\limits^N_{i=1}\omega_{ji}=\sum\limits^N_{i=1}\log\omega_{ji}$,
where $\omega_{ji}>0$, $i=1, \cdots, N$ are determined by
\eqref{omega-i} and \eqref{lambda} with $\mu$ being replaced by
$\mu_j$,  $\sum\limits^N_{i=1}\omega_{ji}=1$,
$\sum\limits^N_{i=1}\omega_{ji}y_i=\mu_j$ for $j=1,2$.

For every $t\in (0,1)$, set
$\omega_{ti}=t\omega_{1i}+(1-t)\omega_{2i}$, $i=1, \cdots, N$.  Then
$\omega_{ti}>0$, $\sum\limits^N_{i=1}\omega_{ti}=1$,
$\sum\limits^N_{i=1}\omega_{ti}y_i=t\mu_1+(1-t)\mu_2\in H$.  Since
$\log x$ is strictly concave in $(0,\infty)$, we have
\[
\log(\omega_{ti})=\log\Big(t\omega_{1i}+(1-t)\omega_{2i}\Big)\ge
t\log\omega_{1i}+(1-t)\log\omega_{2i}~~~i=1,\cdots, N,
\]
and at least one of the inequalities is strict, i.e,
$\log(\omega_{ti})> t\log\omega_{1i}+(1-t)\log\omega_{2i}$ for some
$i$, since $\mu_1\ne\mu_2$ implies $(\omega_{11},\omega_{12},
\cdots, \omega_{1N})\ne (\omega_{21}, \omega_{22},\cdots,
\omega_{2N})$.  Therefore, we get
\[
\sum^N_{i=1}\log(\omega_{ti})>
\sum^N_{i=1}\Big(t\log\omega_{1i}+(1-t)\log\omega_{2i}\Big)=t\log
L(\mu_1)+(1-t)\log L(\mu_2),
\]
which implies
\[
\log L(t\mu_1+(1-t)\mu_2)\ge
\log\prod^N_{i=1}\omega_{ti}=\sum^N_{i=1}\log(\omega_{ti})>t\log
L(\mu_1)+(1-t)\log L(\mu_2),
\]
proving \eqref{concave}.

When $\mu=\bar y$,  an obvious solution to \eqref{lambda} is $\lambda=0$. Since the solution to \eqref{lambda} is unique, we see that $\omega_i=N^{-1}$, and thus, $L(\bar y)=N^{-N}$.  We also notice that
\[
\sup_{\mu}L(\mu)=\sup_{\mu\in H}L(\mu)\le\sup\Big\{\prod^N_{i=1}\omega_i\Big|\omega_i\ge 0, \sum^N_{i=1}\omega_i=1\Big\}=N^{-N}.
\]
The last step is obtained by using the Lagrange multipliers. We omit
the details here. Therefore, we conclude that $L(\bar
y)=\sup\limits_{\mu}L(\mu)=N^{-N}$. \eop

\vspace{10pt}

\noindent{\it Proof of Theorem~\ref{thm1}.} We assume the null hypothesis in \eqref{HH} is true in the proof.

Define $\sigma_n^2=\frac{2(n-2)}{n+1}$ and $S_n^2=\frac{1}{N}\sum\limits_{1\le i<j\le p}\big((n-1)r_{ij}^2-1\big)^2$.  Review that $N=p(p-1)/2$. We have $\sigma_{np}^2=N\sigma_n^2/(n-1)^2$.  Since the distribution of $y_j$'s depends on $n$,
 $\{y_j,\,1\le j\le N\}$ forms an array of random variables.

If the following three conditions are satisfied:  (i). $\frac1{\sigma_n}\max\limits_{1\le j\le N}|y_j-1|=o_p(N^{1/2})$ as $n\to\infty$; (ii).
  $\frac1{N\sigma_n^2}\sum\limits^N_{j=1}(y_j-1)^2\overset{p}{\to} 1$ as $n\to\infty$; (iii).  $\displaystyle{\frac{\sum\limits^N_{j=1}y_j-N}{\sqrt{N\sigma_n^2}}}\overset{d}\to N(0,1)$ as $n\to\infty$,  equivalently, in term of $r_{ij}^2$'s,
\begin{enumerate}
  \item[(C1).] $\displaystyle\frac1{\sigma_n}\max_{1\le i<j\le p}|(n-1)r_{ij}^2-1|=o_p(N^{1/2})$ as $n\to\infty$;
  \item[(C2).] $\displaystyle\frac1{\sigma_n^2}S_n^2\overset{p}{\to} 1$ as $n\to\infty$;
  \item[(C3).] $\displaystyle z_n:=\frac{\sum\limits_{1\le i<j\le p}(n-1)r_{ij}^2-N}{\sqrt{N\sigma_n^2}}\overset{d}\to N(0,1)$ as $n\to\infty$,
\end{enumerate}
we can follow the same procedure as in Owen~\cite{Owen2001} or use Theorem 6.1 in Peng and Schick~\cite{PS2013} to conclude that
\[
\ell(1)=\Big(\frac{\sum\limits_{1\le i<j\le p}(n-1)r_{ij}^2-N}{\sqrt{N\sigma_n^2}}\Big)^2(1+o_p(1))+o_p(1)=z_n^2(1+o_p(1))+o_p(1)
\]
where $\ell(1)$ is defined in \eqref{lmu_0} with $\mu_0=1$. Again, by using condition (C3), we have
\[
\ell_n=z_n^2(1+o_p(1))I(\bar y>0)+o_p(1)=z_n^2(1+o_p(1))I(z_n>0)+o_p(1)\overset{d}\to Z^2I(Z>0)
\]
as $n\to\infty$, where $\ell_n$ is defined in \eqref{l-one}, proving Theorem~\ref{thm1}.

Now we will verify conditions (C1), (C2) and (C3). (C3) has been proved by Chang and Qi~\cite{CQ18}  as we indicate below equation \eqref{schott1}.

Assume $(i,j)$ is a pair of integers with for $1\le i<j\le p$.
It is proved in Schott~\cite{Schott2005} that
\begin{equation}\label{schott2}
E(r_{ij}^2)=\frac{1}{n-1}, ~\text{Var}(r_{ij}^2)=\frac{2(n-2)}{(n+1)(n-1)^2}=\frac{\sigma_n^2}{(n-1)^2},
\end{equation}

From Chang and Qi~\cite{CQ18}, we have
\begin{equation}\label{changqi}
\begin{aligned}
&E(r_{ij}^4)=\frac{3}{(n-1)(n+1)},~E(r_{ij}^6)=\frac{15}{(n-1)(n+1)(n+3)},\\
&E(r_{ij}^8)=\frac{105}{(n-1)(n+1)(n+3)(n+5)}.
\end{aligned}
\end{equation}
By using binomial expansion, we also have
\begin{equation}\label{m4}
m_4:=E\big((r_{ij}^2-\frac{1}{n-1})^4\big)=E(r_{ij}^8)-4\frac{E(r_{ij}^6)}{n-1}+6\frac{E(r_{ij}^4)}{(n-1)^2}-4\frac{E(r_{ij}^2)}{(n-1)^3}+\frac{1}{(n-1)^4}.
\end{equation}

Now we can verify condition (C1). By use of Chebyshev's  inequality, equations \eqref{schott2} and \eqref{changqi}
\begin{eqnarray*}
P(\frac{\max\limits_{1\le i<j\le p}(n-1)r_{ij}^2}{\sigma_n}>\delta N^{1/2})&\le & \sum_{1\le i<j\le p}P(\frac{r_{ij}^2}{\sigma_n}>\frac{\delta N^{1/2}}{n-1})\\
&\le&
\frac{N(n-1)^3}{\delta^4N^{3/2}\sigma_n^3}E(r_{12}^6)\\
&=&O(\frac{1}{N^{1/2}})\to 0
\end{eqnarray*}
as $n\to\infty$ for every $\delta>0$. This implies $\frac{1}{\sigma_n}\max\limits_{1\le i<j\le p}(n-1)r_{ij}^2=o(N^{1/2})$. Hence, we have
\[
\frac1{\sigma_n}\max_{1\le i<j\le p}|(n-1)r_{ij}^2-1|=\frac{1}{\sigma_n}\max_{1\le i<j\le p}(n-1)r_{ij}^2+O(1)= o_p(N^{1/2}),
\]
proving condition (C1).

Below we will use  $(i,j)$ and $(s, t)$ to denote two pair of integers with $1\le i<j\le p$ and $1\le s<t\le p$.
It follows from Theorem 2 in  Veleval and Ignatov~\cite{VI2006} that $\{r_{ij}, ~1\le i<j\le p\}$ are pairwise independent, that is,
 If $(i,j)\ne (s,t)$, then
$r_{ij}$ and $r_{st}$ are independent, thus we have
\[
E\Big(\big((n-1)r_{ij}^2-1\big)^2\big((n-1)r_{st}^2-1\big)^2\Big)=E\Big(\big((n-1)r_{ij}^2-1\big)^2\Big)E\Big(\big((n-1)r_{st}^2-1\big)^2\Big)=
\sigma_n^4.
\]

Since $E(S_n^2)=\sigma_n^2$, we have
\begin{eqnarray*}
E\Big(S_n^2-\sigma^2_n \Big)^2&=&E(S_n^4)-\sigma_n^4\\
&=&\frac{1}{N^2}\sum_{1\le i<j\le p}\sum_{1\le s<t\le p}E\Big(\big((n-1)r_{ij}^2-1\big)^2\big((n-1)r_{st}^2-1\big)^2\Big)-\sigma_n^4.
\end{eqnarray*}
We can classify the summands within the double summation above into
two classes: $N(N-1)$ terms in class 1 when  $(i,j)\ne (s,t)$ and
$N$ terms in class 2 when  $(i,j)=(s,t)$. We see that
\[
E\Big(\big((n-1)r_{ij}^2-1\big)^2\big((n-1)r_{st}^2-1\big)^2\Big)=\sigma_n^4
\]
if $(i,j)\ne (s,t)$ by using the independence, and
\[
E\Big(\big((n-1)r_{ij}^2-1\big)^2\big((n-1)r_{st}^2-1\big)^2\Big)=m_4(n-1)^4
\] if
$(i,j)=(s,t)$, where $m_4$ is given by \eqref{m4}.  Therefore, we
have
\[
E\Big(S_n^2-\sigma^2_n \Big)^2=\frac{1}{N^2}(N(N-1)\sigma_n^4+N m_4(n-1)^4)-\sigma_n^4=\frac{m_4(n-1)^4-\sigma_n^4}{N}.
\]
In view of \eqref{schott2}, \eqref{changqi} and \eqref{m4}, some
tedious calculation shows that
\[
E\Big(S_n^2-\sigma^2_n \Big)^2=\frac{16(n-2)(7n^3-30n^2+11n+60)}{p(p-1)(n+1)^2(n+3)(n+5)}=\frac{4(7n^3-30n^2+11n+60)\sigma_n^4}{p(p-1)(n-2)(n+3)(n+5)},
\]
which implies
\[
E\Big(\frac{S_n^2}{\sigma^2_n}-1\Big)^2=
\frac{4(7n^3-30n^2+11n+60)}{p(p-1)(n-2)(n+3)(n+5)}=O(\frac{1}{p_n^2})\to 0,
\]
as $n\to\infty$, and thus Condition (C2) holds. The proof of
Theorem~\ref{thm1} is complete. \eop

\noindent{\it Proof of Theorem~\ref{thm2}.} Theorem 2 can be proved by using similar arguments in the proof of Theorem~\ref{thm1}. We will continue to use the notation defined in the proof of Theorem~\ref{thm1}.

Under the conditions in Theorem~\ref{thm2}, Chen and Shao~\cite{CS2012} have obtained the following results:
\begin{align}\label{CS-estimate}
    \begin{aligned}
    &E(r_{ij}^2)=\frac{1}{n-1}, ~~E(r_{ij}^4)=\frac{3}{n^2}+O(\frac{1}{n^3}),~~ E(r_{ij}^8)=O(\frac{1}{n^4}); \\
    &E(r_{ij_1}^2r_{ij_2}^2)=\frac{1}{(n-1)^2},~~E(r_{ij_1}^4r_{ij_2}^4)=\frac{9}{n^4}+O(\frac{1}{n^5}) ~~\mbox{ if }j_1\ne j_2,
    \end{aligned}
\end{align}
where $1\le i\ne j\le p$, $1\le i\ne j_1\le p$ and $1\le i\ne j_2\le p$.  It follows from the $C_r$ inequality that
\begin{equation}\label{8th-moment}
E\big((n-1)r_{ij}^2-1\big)^4\le 2^3\big((n-1)^4r_{ij}^8+1\big)=O(1).
\end{equation}

We need to verify conditions (C1), (C2), and (C3) as given in the proof of theorem~\ref{thm1}.   Condition (C1) can be verified similarly by using estimates of moments in \eqref{CS-estimate}, and condition (C3) follows from the central limit theorem
\eqref{schott1}, which is true under the condition of the theorem in virtue of Theorem 2.2 in Chen and Shao~\cite{CS2012}.

Now we proceed to verify condition (C2). First, we have
\begin{equation}\label{2sigma}
\bar{\sigma}_n^2:=E(S_n^2)=E\Big(\big((n-1)r_{12}^2-1\big)^2\Big)=(n-1)^2E\big(r_{12}^4\big)-1=2+O(\frac{1}{n})=\sigma_n^2(1+O(\frac{1}{n}))
\end{equation}
from \eqref{CS-estimate}. Then
\begin{eqnarray*}
E\Big(S_n^2-\bar{\sigma}^2_n \Big)^2&=&E(S_n^4)-\bar{\sigma}_n^4\\
&=&\frac{1}{N^2}\sum_{1\le i<j\le p}\sum_{1\le s<t\le p}E\Big(\big((n-1)r_{ij}^2-1\big)^2\big((n-1)r_{st}^2-1\big)^2\Big)-\bar{\sigma}_n^4.
\end{eqnarray*}
Considering the summands within the double summation above, we see that there are $N(p-2)(p-3)/2$ pairs of sets $\{i,j\}$ and $\{s,t\}$ which are disjoint. For these pairs,
\[
E\Big(\big((n-1)r_{ij}^2-1\big)^2\big((n-1)r_{st}^2-1\big)^2\Big)=\bar{\sigma}_n^4.
\]
because $r_{ij}$ and $r_{st}$ are independent. For all other $N^2-N(p-2)(p-3)/2$ pairs, corresponding summands are dominated by
\begin{eqnarray*}
E\Big(\big((n-1)r_{ij}^2-1\big)^2\big((n-1)r_{st}^2-1\big)^2
&\le&\sqrt{E\big((n-1)r_{ij}^2-1\big)^4
E\big((n-1)r_{st}^2-1\big)^4}\\
&=&E\big((n-1)r_{ij}^2-1\big)^4\\
&=&O(1)
\end{eqnarray*}
from the Cauchy-Schwarz inequality and equation \eqref{8th-moment}.  Therefore, we have
\begin{eqnarray*}
E\Big(S_n^2-\bar{\sigma}^2_n \Big)^2&=& \frac{1}{N^2}\Big(\frac12{N(p-2)(p-3)}\bar{\sigma}_n^4+
O\big(N^2-\frac12N(p-2)(p-3)\big)\Big)-\bar{\sigma}_n^4\\
&=&O\Big(\frac1{N^2}\big(N^2-\frac12N(p-2)(p-3)\big)\Big)\\
&=&O(\frac1p)\to 0
\end{eqnarray*}
as $n\to\infty$. In the estimation above, we also use the fact that $\bar{\sigma}_n^2\sim \sigma_n^2\to 2$ from \eqref{2sigma}. Therefore, we have
\[
E\Big(\frac{S_n^2}{\bar{\sigma}^2_n}-1\Big)^2\to 0
\]
as $n\to\infty$. This yields $\displaystyle{\frac{S_n^2}{\bar{\sigma}^2_n}\overset{p}\to 1}$, which together with \eqref{2sigma} implies condition (C2).  This completes the proof of the theorem.   \eop

\noindent{\it Proof of Equation \eqref{1-in-probability}}.
In the proofs of Theorems~\ref{thm2} and \ref{thm2}, we have obtained that $S_n/\sigma_n^2\overset{p}{\to} 1$, which implies $S_n\overset{p}{\to} 2$ since $\sigma_n^2\to 2$.
From condition (C3), we have
\[
\frac{1}{N}\big(\sum_{1\le i<j\le p}(n-1)^2r_{ij}^2-N\big)
\overset{p}{\to} 0.
\]
Therefore, we get
\begin{eqnarray*}
 \frac{2(n-1)(n+1)}{3(p-1)(p+4)}\sum_{1\le i<j\le p}r_{ij}^4&=&\frac{1+o(1)}{3N}\sum_{1\le i<j\le p}(n-1)^2r_{ij}^4\\
 &=& \frac{1+o(1)}{3}\Big(S_n^2+1+\frac{2}{N}\big(\sum_{1\le i<j\le p}(n-1)^2r_{ij}^2-N\big)\Big)\\
 &\overset{p}\to &1
\end{eqnarray*}
as $n\to\infty$, proving \eqref{1-in-probability}.  \eop

\section*{Acknowledgements}

The authors would like to thank the Editor, Associate Editor, and referee for reviewing the manuscript and providing valuable comments.
The research of Yongcheng Qi was supported in part by NSF Grant DMS-1916014.

\section*{Disclosure statement}

No potential conflict of interest was reported by the authors.

\newpage

\begin{table}[H]
\caption{Sizes and Powers of tests,   Normal distribution} \label{table-norm} \centering
{\footnotesize
\begin{tabular}{|ll|ccc|ccc|ccc|}
 \hline
 & &\multicolumn{3}{c|}{size ($\rho=0$)} &\multicolumn{3}{c|}{power ($\rho=0.02$)}  & \multicolumn{3}{c|}{power ($\rho=0.05$)}\\
   \cline{3-11}
 $n$ & $p$& $\bar{\ell}_n$ & $t_{np}^c$ & $t_{np}^*$ & $\bar{\ell}_n$ & $t_{np}^c$ & $t_{np}^*$& $\bar{\ell}_n$ & $t_{np}^c$ & $t_{np}^*$\\
   \hline
 20 & 10  & 0.0593 & 0.0486 & 0.0598 & 0.0646 & 0.0551 & 0.0657 & 0.0974 & 0.0863 & 0.0987 \\
    & 20  & 0.0609 & 0.0565 & 0.0618 & 0.0670 & 0.0626 & 0.0675 & 0.1357 & 0.1281 & 0.1366 \\
    & 50  & 0.0594 & 0.0584 & 0.0596 & 0.0836 & 0.0810 & 0.0842 & 0.3079 & 0.3030 & 0.3096 \\
    & 100 & 0.0524 & 0.0519 & 0.0525 & 0.1156 & 0.1142 & 0.1159 & 0.5584 & 0.5571 & 0.5588 \\
    \hline
 50 & 10  & 0.0596 & 0.0484 & 0.0603 & 0.0750  & 0.064 & 0.0744 & 0.1660  & 0.1478 & 0.1656 \\
    & 20  & 0.0533 & 0.0481 & 0.0525 & 0.0855 & 0.0790 & 0.0859 & 0.3257 & 0.3098 & 0.3244 \\
    & 50  & 0.0508 & 0.0489 & 0.0508 & 0.1457 & 0.1410 & 0.1461 & 0.7372 & 0.7320  & 0.7362 \\
    & 100 & 0.0536 & 0.0519 & 0.0535 & 0.2684 & 0.2660 & 0.2684 & 0.9609 & 0.9604 & 0.9608 \\
    \hline
 100 & 10 & 0.0608 & 0.0504 & 0.0607 & 0.0894 & 0.0743 & 0.0882 & 0.3305 & 0.3025 & 0.3285 \\
    & 20  & 0.0581 & 0.0511 & 0.0573 & 0.1265 & 0.1167 & 0.1260 & 0.6520 & 0.6345 & 0.6495 \\
    & 50  & 0.0523 & 0.0494 & 0.0519 & 0.2735 & 0.2675 & 0.2723 & 0.9816 & 0.9809 & 0.9816 \\
    & 100 & 0.0516 & 0.0506 & 0.0515 & 0.5711 & 0.5669 & 0.5703 & 0.9997 & 0.9997 & 0.9997 \\
    \hline
\end{tabular}
}
\end{table}

\begin{table}[H]
\caption{Sizes and Powers of tests, Uniform distribution} \label{table-unif} \centering
{\footnotesize
\begin{tabular}{|ll|ccc|ccc|ccc|}
 \hline
 & &\multicolumn{3}{c|}{size ($\rho=0$)} &\multicolumn{3}{c|}{power ($\rho=0.02$)}  & \multicolumn{3}{c|}{power ($\rho=0.05$)}\\
   \cline{3-11}
 $n$ & $p$& $\bar{\ell}_n$ & $t_{np}^c$ & $t_{np}^*$ & $\bar{\ell}_n$ & $t_{np}^c$ & $t_{np}^*$& $\bar{\ell}_n$ & $t_{np}^c$ & $t_{np}^*$\\
   \hline
 20 & 10 & 0.0610           & 0.0529     & 0.0615     & 0.0607         & 0.0525     & 0.0626     & 0.0913         & 0.0802     & 0.0934     \\
    & 20  & 0.0600           & 0.0571     & 0.0615     & 0.0668         & 0.0621     & 0.0686     & 0.1230          & 0.1147     & 0.1243     \\
    & 50  & 0.0595         & 0.0582     & 0.0607     & 0.0806         & 0.0785     & 0.0812     & 0.2680          & 0.2647     & 0.2693     \\
    & 100 & 0.0585         & 0.0582     & 0.0598     & 0.1099         & 0.1093     & 0.1109     & 0.5183         & 0.5171     & 0.5202     \\
    \hline
 50  & 10  & 0.0572         & 0.0489     & 0.0575     & 0.0773         & 0.0655     & 0.0774     & 0.1744         & 0.1536     & 0.1749     \\
    & 20  & 0.0606         & 0.0548     & 0.0600       & 0.0809         & 0.0734     & 0.0816     & 0.3061         & 0.2920      & 0.3062     \\
    & 50  & 0.0559         & 0.0530      & 0.0562     & 0.1358         & 0.1319     & 0.1359     & 0.7449         & 0.7387     & 0.7449     \\
    & 100 & 0.0526         & 0.0518     & 0.0527     & 0.2417         & 0.2378     & 0.2422     & 0.9741         & 0.9736     & 0.9741     \\
    \hline
 100 & 10  & 0.0596         & 0.0497     & 0.0583     & 0.0842         & 0.0710      & 0.0848     & 0.3225         & 0.2934     & 0.3200       \\
    & 20  & 0.0589         & 0.0528     & 0.0580      & 0.1269         & 0.1152     & 0.1271     & 0.6430          & 0.6265     & 0.6397     \\
    & 50  & 0.0516         & 0.0489     & 0.0522     & 0.2650          & 0.2587     & 0.2646     & 0.9889         & 0.9880      & 0.9887     \\
    & 100 & 0.0504         & 0.0494     & 0.0504     & 0.5567         & 0.5525     & 0.5566     & 1.0000              & 1.0000          & 1.0000 \\
    \hline
\end{tabular}
}
\end{table}

\begin{table}[H]
\caption{Sizes and Powers of tests, Exponential distribution} \label{table-exp} \centering
{\footnotesize
\begin{tabular}{|ll|ccc|ccc|ccc|}
 \hline
 & &\multicolumn{3}{c|}{size ($\rho=0$)} &\multicolumn{3}{c|}{power ($\rho=0.02$)}  & \multicolumn{3}{c|}{power ($\rho=0.05$)}\\
   \cline{3-11}
 $n$ & $p$& $\bar{\ell}_n$ & $t_{np}^c$ & $t_{np}^*$ & $\bar{\ell}_n$ & $t_{np}^c$ & $t_{np}^*$& $\bar{\ell}_n$ & $t_{np}^c$ & $t_{np}^*$\\
   \hline
 20  & 10  & 0.0718         & 0.0636     & 0.0745     & 0.0799         & 0.0734     & 0.0845     & 0.1419         & 0.1338     & 0.1486     \\
    & 20  & 0.0638         & 0.0610      & 0.0670      & 0.0867         & 0.0859     & 0.0924     & 0.2135         & 0.2116     & 0.2221     \\
    & 50  & 0.0581         & 0.0586     & 0.0612     & 0.1224         & 0.1266     & 0.1301     & 0.3982         & 0.4040      & 0.4080      \\
    & 100 & 0.0602         & 0.0622     & 0.0632     & 0.1901         & 0.1955     & 0.1975     & 0.5774         & 0.5848     & 0.5871     \\
    \hline
50  & 10  & 0.0734         & 0.0620      & 0.0751     & 0.0950          & 0.0840      & 0.0975     & 0.2288         & 0.2156     & 0.2335     \\
    & 20  & 0.0650          & 0.0606     & 0.0660      & 0.1158         & 0.1114     & 0.1205     & 0.3834         & 0.3800       & 0.3898     \\
    & 50  & 0.0612         & 0.0608     & 0.0643     & 0.1877         & 0.1912     & 0.1953     & 0.7102         & 0.7140      & 0.7175     \\
    & 100 & 0.0589         & 0.0607     & 0.0618     & 0.3230          & 0.3289     & 0.3310      & 0.9079         & 0.9107     & 0.9111     \\
    \hline
100 & 10  & 0.0701         & 0.0604     & 0.0698     & 0.1048         & 0.0936     & 0.1088     & 0.3636         & 0.3428     & 0.3689     \\
    & 20  & 0.0646         & 0.0595     & 0.0647     & 0.1555         & 0.1492     & 0.1578     & 0.6350          & 0.6281     & 0.6400       \\
    & 50  & 0.0636         & 0.0622     & 0.0653     & 0.3108         & 0.3125     & 0.3187     & 0.9481         & 0.9486     & 0.9495     \\
    & 100 & 0.0574         & 0.0583     & 0.0600       & 0.5823         & 0.5878     & 0.5899     & 0.9965         & 0.9966     & 0.9966    \\
    \hline
\end{tabular}
}
\end{table}

\begin{table}[H]
\caption{Sizes and Powers of tests, Mixture distribution} \label{table-mix} \centering
{\footnotesize
\begin{tabular}{|ll|ccc|ccc|ccc|}
 \hline
 & &\multicolumn{3}{c|}{size ($\rho=0$)} &\multicolumn{3}{c|}{power ($\rho=0.02$)}  & \multicolumn{3}{c|}{power ($\rho=0.05$)}\\
   \cline{3-11}
 $n$ & $p$& $\bar{\ell}_n$ & $t_{np}^c$ & $t_{np}^*$ & $\bar{\ell}_n$ & $t_{np}^c$ & $t_{np}^*$& $\bar{\ell}_n$ & $t_{np}^c$ & $t_{np}^*$\\
   \hline
 20  & 10  & 0.0650          & 0.0547     & 0.0664     & 0.0682         & 0.0594     & 0.0680      & 0.1042         & 0.0921     & 0.1052     \\
    & 20  & 0.0618         & 0.0565     & 0.0630      & 0.0704         & 0.0650      & 0.0706     & 0.1561         & 0.1477     & 0.1572     \\
    & 50  & 0.0605         & 0.0584     & 0.0608     & 0.0974         & 0.0955     & 0.0983     & 0.3381         & 0.3353     & 0.3394     \\
    & 100 & 0.0688         & 0.0677     & 0.0692     & 0.1395         & 0.1381     & 0.1407     & 0.5979         & 0.5974     & 0.5998     \\
    \hline
50  & 10  & 0.0624         & 0.0520      & 0.0616     & 0.0778         & 0.0657     & 0.0781     & 0.1991         & 0.1790      & 0.1988     \\
    & 20  & 0.0602         & 0.0548     & 0.0606     & 0.0958         & 0.0869     & 0.0952     & 0.3832         & 0.3713     & 0.3828     \\
    & 50  & 0.0587         & 0.0559     & 0.0578     & 0.1643         & 0.1595     & 0.1638     & 0.7827         & 0.7792     & 0.7824     \\
    & 100 & 0.0569         & 0.0561     & 0.057      & 0.3086         & 0.3063     & 0.3097     & 0.9669         & 0.9666     & 0.9672     \\
    \hline
100 & 10  & 0.0610          & 0.0510      & 0.0615     & 0.0990          & 0.0845     & 0.0961     & 0.3777         & 0.3484     & 0.3745     \\
    & 20  & 0.0599         & 0.0520      & 0.0592     & 0.1395         & 0.1299     & 0.1380      & 0.7200           & 0.7055     & 0.7183     \\
    & 50  & 0.0556         & 0.0522     & 0.0555     & 0.3137         & 0.3080      & 0.3137     & 0.9888         & 0.9887     & 0.9888     \\
    & 100 & 0.0561         & 0.0545     & 0.0559     & 0.6389         & 0.6359     & 0.6400       & 1.0000              & 1.0000          & 1.0000 \\
    \hline
\end{tabular}
}
\end{table}

\begin{table}[H]
\caption{Sizes and Powers of tests, Sum} \label{table-sum} \centering
{\footnotesize
\begin{tabular}{|ll|ccc|ccc|ccc|}
 \hline
 & &\multicolumn{3}{c|}{size ($\rho=0$)} &\multicolumn{3}{c|}{power ($\rho=0.02$)}  & \multicolumn{3}{c|}{power ($\rho=0.05$)}\\
   \cline{3-11}
 $n$ & $p$& $\bar{\ell}_n$ & $t_{np}^c$ & $t_{np}^*$ & $\bar{\ell}_n$ & $t_{np}^c$ & $t_{np}^*$& $\bar{\ell}_n$ & $t_{np}^c$ & $t_{np}^*$\\
   \hline
 20  & 10  & 0.0601         & 0.0520      & 0.0607     & 0.0625         & 0.0552     & 0.0634     & 0.0997         & 0.0884     & 0.1019     \\
    & 20  & 0.0556         & 0.0499     & 0.0558     & 0.0635         & 0.0584     & 0.0637     & 0.1402         & 0.1303     & 0.1401     \\
    & 50  & 0.0520          & 0.0499     & 0.0520      & 0.0804         & 0.0773     & 0.0807     & 0.3066         & 0.3031     & 0.3072     \\
    & 100 & 0.0592         & 0.0583     & 0.0593     & 0.1172         & 0.1153     & 0.1173     & 0.5509         & 0.5484     & 0.5515     \\
    \hline
 50  & 10  & 0.0587         & 0.0493     & 0.0586     & 0.0723         & 0.0610      & 0.0723     & 0.1702         & 0.1501     & 0.1684     \\
    & 20  & 0.0558         & 0.0497     & 0.0557     & 0.0843         & 0.0762     & 0.0838     & 0.3315         & 0.3144     & 0.3290      \\
    & 50  & 0.0537         & 0.0507     & 0.0541     & 0.1429         & 0.1385     & 0.1429     & 0.7446         & 0.7377     & 0.7433     \\
    & 100 & 0.0557         & 0.0540      & 0.0555     & 0.2637         & 0.2614     & 0.2642     & 0.9550          & 0.9545     & 0.9550      \\
    \hline
 100 & 10  & 0.0635         & 0.0552     & 0.0640      & 0.0897         & 0.0758     & 0.0892     & 0.3237         & 0.2955     & 0.3187     \\
    & 20  & 0.0593         & 0.0527     & 0.0586     & 0.1211         & 0.1097     & 0.1197     & 0.6540          & 0.6360      & 0.6513     \\
    & 50  & 0.0480          & 0.0459     & 0.0479     & 0.2716         & 0.2644     & 0.2706     & 0.9803         & 0.9791     & 0.9801     \\
    & 100 & 0.0551         & 0.0537     & 0.0549     & 0.5719         & 0.5690      & 0.5713     & 1.0000              & 1.0000          & 1.0000  \\
    \hline
\end{tabular}
}
\end{table}

\end{document}